\documentclass[11pt]{article}
\usepackage{amssymb,amsmath,amsfonts}
\usepackage[active]{srcltx}
\usepackage{color}
\usepackage{ulem}

\usepackage{epsfig}
\newtheorem{Lemma}{Lemma}[section]
\newtheorem{Theorem}[Lemma]{Theorem}

\newtheorem{Corollary}[Lemma]{Corollary}
\newtheorem{Remark}[Lemma]{Remark}

\newtheorem{Definition}[Lemma]{Definition}

\def\nc{\newcommand}
\def\rnc{\renewcommand}

\newcommand{\qed}{\nobreak \ifvmode \relax \else
      \ifdim\lastskip<1.5em \hskip-\lastskip
      \hskip1.5em plus0em minus0.5em \fi \nobreak
      \vrule height0.75em width0.5em depth0.25em\fi}
\nc{\Halmos}    {\ \raisebox{0.6ex}  {\framebox[0.9ex]{
			   \rule[0ex]{0ex}{0.5ex} }}}

\def\theequation{\arabic{section}.\arabic{equation}}

\nc{\mylabel}[1]{\label{#1}}
\nc{\myref}[1]{{\mbox{\bf [#1]}} \ref{#1}}
\nc{\isep}{\setlength{\itemsep}{-1mm}}
\nc{\ds}{\displaystyle}

\nc{\markred}{\textcolor{red}}
\nc{\markgreen}{\textcolor{green}}

\nc{\mat}[2]  {\left(  \! \begin{array}{#1} #2 \end{array}\! \right)}
\nc{\matabcd}{\mat{cc}{a&b\\c&d}}

\nc{\rr}    {\rightarrow}

  \nc{\tf}{{\tilde{f}}}
  \nc{\tg}{{\tilde{g}}}
\nc{\tx}{{\tilde{x}}} \nc{\ty}{{\tilde{y}}}
\nc{\tF}{{\tilde{F}}}

\newlength{\hhgt}
\nc{\mybar}[1]{
    \settoheight{\hhgt}{{#1}}
    \addtolength{\hhgt}{0.2em}
    \overline{\rule[0ex]{0ex}{\hhgt}{}{#1}}}

\nc{\rstrut} {{\rule[-1ex]{0ex}{2.2ex}{}}}

\nc{\st} {{\rule[-1ex]{0ex}{3ex}{}}}
\nc{\hugestrut} {{\rule[-0.5ex]{0ex}{3.5ex}{}}}
\nc{\mystrut} {{\rule[0ex]{0ex}{2.5ex}{}}}
\nc{\mstrut} {{\rule[0ex]{0ex}{2.1ex}{}}}
\nc{\nstrut} {{\rule[0ex]{0ex}{1.7ex}{}}}
\nc{\qstrut} {{\rule[0ex]{0ex}{1.8ex}{}}}
\nc{\pstrut} {{\rule[0ex]{0ex}{1.5ex}{}}}

\nc{\etaD}{\eta_{\qstrut D}}
\nc{\etasp}{\eta_{\rm sp}}
\nc{\rsp}{r_{\rm sp}}

\nc{\CREP}{{g}}
\rnc{\Re}    {{\rm Re\,}}
\rnc{\Im}    {{\rm Im\,}}
\nc{\re}{{\rm Re}} \nc{\im}{{\rm Im}} 

\nc{\dInner}{ d^{\mstrut \, \rm in}}
\nc{\dRiem}{ d^{\rm Rie}}

\nc{\rank}   {{\rm rank}}
\nc{\ess}    {{\rm ess}}
\nc{\osc}    {{\rm osc}}
\nc{\esssup} {{\rm ess\, sup}}

\nc{\Ra}{\Rightarrow}
\nc{\hra}{\hookrightarrow}
 \nc{\Arg}   {{\rm Arg}}
 \nc{\Aper}   {{\rm Aper}}

\nc{\const}{{\rm const}} \nc{\diam}{{\rm diam}} 

\nc{\bfone}{{\bf 1}}  
\nc{\Id}{{\it Id}}  

\nc{\bfM}  {{\bf M}}
\nc{\bfMo}  {{\bf M}_{\; \omega}}
\nc{\bfMno}  {{\bf M}_{\; \omega}^{(n)}}

\nc{\bfh}  {{\bf h}}
\nc{\bfhs}  {{\bf h}^*}
\nc{\bfho}  {{\bf h}_\omega}
\nc{\bfhto}  {{\bf h}_{\tau \omega}}

\nc{\bfp}  {{\bf p}}
\nc{\bfpo}  {{\bf p}_\omega}

\nc{\bfpi}  { \mbox{\boldmath {$\pi$}}}

\nc{\bfphi}  {{\bf \phi}}

\nc{\LRA}{\Leftrightarrow}

 \nc{\Int}{{\rm Int\;}} \nc{\Imm}{\mbox{Im}}
\nc{\half}{\frac{1}{2}} \nc{\DDD}{D}

\nc{\calP}{{\cal P}}
\nc{\calA}{{\cal A}}
\nc{\calB}{{\cal B}}
\nc{\calE}{{\cal E}}
\nc{\calC}{{\cal C}}
\nc{\calK}{{\cal K}}
\nc{\calL}{{\cal L}}
\nc{\calR}{{\cal R}}
\nc{\calM}{{\cal M}}
\nc{\calT}{{\cal T}}
\nc{\ellone}{{\{\ell=1\}}}
\nc{\Cellone}{{{\C} \cap \ellone}}
\nc{\D}{{\cal D}} 

\nc{\oR}{\overline{R}}
\nc{\uR}{\underline{R}}
\nc{\calCS}      {\calC^{\pstrut *}}
\nc{\cone}{\,\Rplus\,}
\nc{\C}       {{K}}   
\nc{\CS}      {\C^{\pstrut *}}
\nc{\CP}      {\C^{\pstrut '}}
\nc{\XS}      {X^{\pstrut *}}
\nc{\KS}      {K^{\pstrut *}}

\nc{\RS}      {\RR\stt{2}^{\nstrut *}}

\nc{\kk}{{\mathbb K}}   
\nc{\CCP} {{\mathbb C}{\rm P}}
\nc{\RRP} {{\mathbb R}{\rm P}}
\nc{\CC}{{\mathbb C}} 
\nc{\CCS}{{\mathbb C}^*} 
\nc{\HH}{{\mathbb H}} 
\nc{\CCbf}{{\mathbb C}} 
\nc{\DD}{{\mathbb D}} \nc{\EEE}{\mbox{$\mathbb E$}}
\nc{\RR}{{{\mathbb R}}}
\nc{\RRs}{{\mathbb R}} \nc{\NN}{{\mathbb N}} \nc{\NNs}{{\mathbb N}}
\nc{\ZZ}{\mbox{$\mathbb Z$}} \nc{\AutDD}    {{\rm Aut(\DD;\diag \DD)}}
\nc{\hatH}{\widehat{H}}
\nc{\hatRR}{\widehat{\RR}}
\nc{\hatRn}{\widehat{\RR}^n}
\nc{\hatCC}{\widehat{\CC}}
\nc{\Rplus}{\RR_+}
\nc{\RRp}{\RR_+}

\nc{\stt}[1]{{\rule[0ex]{0ex}{#1ex}{}}}

\newcommand{\dual}[1]{\la #1 \ra}

\nc{\la}{\langle}
\nc{\ra}{\rangle}
\nc{\barA}{\overline{A}} 
\nc{\barM}{\overline{M}} 
\nc{\bara}{\overline{a}} 
\nc{\barc}{\overline{c}} 
\nc{\bard}{\overline{d}} 
\nc{\barz}{\overline{z}} 
\nc{\baru}{\overline{u}} 
\nc{\barb}{\overline{b}} \nc{\barDD}{{\overline{\DD}}}
\nc{\barK}{{\overline{K}}} \nc{\bRs}{\overline{\RR}_+((s))}
\nc{\UbarU}{U\times \overline{U}} \nc{\XbarX}{X\times \overline{X}}
\nc{\barmu}{\overline{\mu}}
\nc{\barla}{\overline{\lambda}}

\nc{\Mn}     {M^{(n)}}
\nc{\pin}     {\pi^{(n)}}
\nc{\lMn}[1] {\la \ell, M^{(n)} #1 \ra}
\nc{\mMn}[1] {\la m, M^{(n)} #1 \ra}

\nc{\hatOmega}{\widehat{\Omega}}
\nc{\hatU}{\widehat{U}}
\nc{\hatA}{\widehat{A}}
\nc{\hatB}{\widehat{B}}
\nc{\hatP}{\widehat{P}}
\nc{\hatQ}{\widehat{Q}}
\nc{\hatE}{\widehat{E}}
\nc{\hatF}{\widehat{F}}
\nc{\hatD}{\widehat{D}}
\nc{\hatx}{\widehat{x}}
\nc{\haty}{\widehat{y}}
\nc{\hatd}{\widehat{d}} \nc{\hatg}{g_\htinyD} \nc{\hatv}{\widehat{v}}
\nc{\Dhatf}{D\widehat{f}} \nc{\Dpsihat}{D\widehat{\psi}_t}
\nc{\hatf}{{\widehat{f}}}

\nc{\Cl}{{\rm Cl\;}} \nc{\len}{\mbox{len\,}} \nc{\diag}{\mbox{diag\ }}
\nc{\diagK}{\mbox{diag}(K)}

\nc{\psec}{\partial_{\Sigma}}
\nc{\Intsec}{{\rm Int}_{\Sigma}}

\nc{\longr}{\longrightarrow}

\begin{document}
\title{ A uniform contraction principle for bounded
Apollonian embeddings.
}
\author{Lo\"ic Dubois and Hans Henrik Rugh.\\
    Helsinki University\footnote{This research was partially 
    funded by the European Research Council.},
    Finland.
    University of Cergy-Pontoise,
    CNRS UMR 8088, France. 
    }
\date {\today}
 \maketitle
\begin{abstract}
Let $\hatH=H \cup \{\infty\}$ denote the standard one-point completion of
a real Hilbert space $H$. 
Given any non-trivial proper sub-set $U\subset \hatH$ one may
define the so-called `Apollonian' metric $d_U$ on $U$.
When $U\subset V \subset \hatH$  are nested proper subsets
we show that their associated Apollonian metrics satisfy
the following 
uniform contraction principle: Let $\Delta=\diam_{V}(U) \in [0,+\infty]$ be
the diameter of the smaller subsets with respect to the large.
Then  for every $x,y\in U$ we have
  \[  d_V(x,y) \leq \tanh \frac{\Delta}{4} \ \ d_U(x,y) .\]
In dimension one, this contraction principle was established by
 Birkhoff \cite{Bir57} for
the Hilbert metric of finite segments on $\RRP^1$.
  In dimension two it was shown by Dubois in \cite{Dub09}
  for subsets of the Riemann sphere
  $\hatCC\sim\widehat{\RR^2}$.
  It is new in the generality stated here.
\end{abstract}

\section{Introduction and results}
There are striking similarities between the projective group for
the real or complex projective lines and
the conformal group of the one-point completion of a real Hilbert space of 
dimension at least 3.
In the first case, the group consists
of M\"obius maps of the form $z\mapsto \frac{az+b}{cz+d}$
and in the second
it is generated by linear isometries, homotheties and the inversion,
 corresponding to M\"obius transformations
supplemented with a complex conjugation.
In both cases one
 needs at least 4 points to define
a group invariant quantity, i.e.\ the cross-ratio. 
Fixing a subset $U$ whose complement contains at least 2 points,
the logarithm
of cross-ratios may then be used
to construct a (semi-)metric on
$U$. 
On the interval $I=(-1,1)$, there is a unique (up to a constant) distance invariant under
M\"obius transformations preserving $I$. This is precisely the restriction of the Poincar\'e metric
$2|dz|(1-|z|^2)^{-1}$ on the unit disk in the complex plane.
In the case of Hilbert spaces of higher dimensions
one may derive the so-called `Apollonian metric' (see below).
This latter metric was first introduced for 
$\hatRn=\RR^n\cup \{\infty\}$  by Barbilian \cite{Bar34}
and later rediscovered by Beardon \cite{Bea98}.

   From a dynamical point of view it is of interest to know how 
a subset $U$ metrically embed into a larger subset
$V$ with respect to the associated metrics $d_U$ and $d_V$
(see below for more precise statements).
It is straight-forward from definitions that the injection
$i:(U,d_U) \hookrightarrow (V,d_V)$ is non-expanding. More interesting, however,
is that it verifies a very general uniform contraction principle
(UCP):
If $\Delta=\diam_V(U)<+\infty$, i.e.\  the embedding of $U$ has bounded
diameter in the larger domain $V$, then the injection is a strict
contraction with a Lipschitz constant bounded by $\tanh \frac{\Delta}{4}$.
This is the same formula which appears in Birkhoff's work on the Hilbert
metric \cite{Bir57}.
We give below the (surprisingly simple) proof of the UCP
for the general case and in section
	\ref{sec applic} some simple
dynamical systems applications.

There is no particular reasons for sticking to finite dimension,
so in the following let $H$ be any real Hilbert space.
We write 
 $\dual{\cdot,\cdot}$ for the scalar product and $\|\cdot\|$ for the
 norm on $H$. 
Let $\hatH=H\cup \{\infty\}$ be a one point completion of $H$ in which
the open sets containing $\infty$ are of the
form $\{\infty\} \cup F^c$ with $F$ a bounded closed set.
With this convention
$\hatH$ is compact iff $H$ is finite dimensional.
The space $(\hatH,\hatd)$ is a complete metric space of diameter one
with respect to the  metric:
\begin{equation}
 \hatd(x,y) = 
      \frac{\|x-y\|}{\sqrt{1+\dual{x,x}}
      \sqrt{1+\dual{y,y}}} , \ \ \
  \hatd(\infty,y) = \frac{1} {\sqrt{1+\dual{y,y}}} \ .
\end{equation}

\begin{Definition}
   \mylabel{cross}
   Given four points $\{x_1,x_2,u_1,u_2\} \in \hatH$ such that
   $\{x_1,x_2\}$ and $\{u_1,u_2\}$ are disjoint we define their
   cross-ratio to be~:
   \begin{equation}
      [x_1,x_2;u_1,u_2] \equiv  \frac
            {\|x_2-u_1\| \; \|x_1-u_2\|}
            {\|x_1-u_1\| \; \|x_2-u_2\|}.
   \end{equation}
   Here, $\|\cdot\|$ denotes the Hilbert norm in $H$ and we adapt
   usual conventions for dealing with the point at $\infty$.
   When $U\subset \hatH$ is a proper subset (by proper we
   mean that $U$ and $U^c$ are both
   non-empty) one defines 
   the Apollonian (semi,pseudo-)distance between points $x_1,x_2\in U$:
   \begin{equation}
                d_{U}(x_1,x_2) = \sup_{u_1,u_2\in U^c}
		   \log [x_1,x_2;u_1,u_2] \in [0,+\infty]
   \end{equation}
\end{Definition}
We denote by $GM(\hatH)$ the general conformal group which acts
continuously upon $(\hatH,\hatd)$ and is
generated by the set of isometries, homotheties
(both fixing $\infty$) and the inversion
(which exchanges the origin
and $\infty$):
\begin{equation}
    I(x)= \frac{x}{\dual{x,x}},  \ \ \
    I(0)=\infty \ \ \ \mbox{and} \ \ \
    I(\infty)=0. \mylabel{inversion}
\end{equation}
When $\dim H \geq 3$ 
the Liouville theorem (see e.g.\ \cite{Nev60}) shows that
any conformal map is in $GM(\hatH)$.
In dimension 1 or 2, it is the M\"obius group (supplemented with
complex conjugation in the 2 dimensional case).
That $d_U$ is $GM(\hatH)$ invariant 
is trivial for isometries and homotheties
and in the case of inversions it
follows easily from the formula
$\ds \|I(x)-I(y)\|=\frac{\|x-y\|}{\|x\| \; \|y\|}$
 (with some care taken with
 respect to the point at infinity).  From the cross-ratio identity
$[x,z;u,v]=[x,y;u,v] \ [y,z;u,v]$ and taking sup
in the right order one also sees that $d_U$ verifies the triangular
inequality. When $U^c$ has non-empty interior 
$d_U$ is a genuine metric, but
in the general case it need not distinguish points.
We refer to e.g. \cite[Chapter 3]{Bea98} and \cite{Has04}
 for further details on the geometry of this metric.
Our main result is the following:

\begin{Theorem} 
    \mylabel{thm contract} [Main Theorem]
    Let $U \subset V \subset \hatH$ 
     be non-empty proper subsets
   with 
 $d_U$ and $d_V$ being the associated Apollonian metrics.
Let $ \Delta =  \sup_{u_1,u_2\in U} d_V(u_1,u_2)$
be the diameter of 
the smaller subset within the larger.
Then for every $x_1,x_2\in U$:
\begin{equation}
      d_V(x_1,x_2) \leq \left( \tanh \frac{\Delta}{4} \right) \ d_U(x_1,x_2).
\end{equation}
If $\diam_V(U)<+\infty$, the embedding
$i: (U,d_V) \hookrightarrow (U,d_U)$ is 
a uniform contraction.\\
\end{Theorem}

Proof:
We will base our proof 
upon Birkhoff's inequality \cite{Bir57} for cross-ratios
on the projective real line. It is, in fact, a special case
of our main theorem when $n=1$.
We will use it in the
following version:
Let $K=(a_1,a_2)$ be a non-empty open sub-interval
of $J=(0,+\infty)$ .
The Hilbert distance of $s_1,s_2\in K$ relative to $K$ and $J$ are
given by:
\[ 
d_{K}(s_1,s_2) = \left|\mstrut \log [s_1,s_2;a_1,a_2] \right| 
\ \ \mbox{and} \ \ 
d_{J}(s_1,s_2) = \left|\mstrut  \log  \frac{s_2}{s_1} \right|
.\]
The quantity
$\ds \Delta= \diam_{J}(K) = 
\log \frac{a_2}{a_1}  \in (0,+\infty]$ 
measures
the diameter of $K$ for the $J$-metric.  Birkhoff \cite[p.220]{Bir57} 
showed the fundamental inequality :
\begin{equation}
d_{J}(s_1,s_2) \leq \left( \tanh \frac{\Delta}{4} \right) \ d_{K}(s_1,s_2),
   \ \ \ \forall \  s_1,s_2\in K. 
\label{Birkhoff ineq}
\end{equation}
Proof of (\ref{Birkhoff ineq}): 
 It suffices to show this for $s_1$ and $s_2$ infinitesimally close.
 So we differentiate with respect to $s_2$ at 
 $s_2=s_1=s\in (a_1,a_2)$ and search for the
 optimal value of $\theta>0$ so that for every $a_1<s<a_2$:
 $\ds \frac{1}{s} \leq \theta \frac{a_1-a_2}{(s-a_1)(a_2-s)}$,
 or equivalently 
  \begin{equation}
  \theta \  \geq \  \inf_{a_1<s<a_2} \frac{(s-a_1)(a_2-s)}{s(a_2-a_1)}.
  \end{equation}
 The minimum value is at $\ds s=\sqrt{a_1a_2}$ and equals
 $\ds \theta_{\rm min} = \frac{\sqrt{a_2}-\sqrt{a_1}}{\sqrt{a_2}+\sqrt{a_1}}=
 \tanh \frac{\log (a_2/a_1)}{4}$ which is therefore
 the desired contraction constant.

 Now, returning to the general case 
 let  $x_1,x_2\in U$ be distinct points.
 We have $d_V(x_1,x_2) \leq d_U(x_1,x_2)$ since
 the sup in the latter case is over a larger set.
So we may assume that $\Delta=\diam_V(U)<+\infty$
and also that $0<d_V(x_1,x_2)\leq d_U(x_1,x_2)<+\infty$
(or else the statement is trivial).
 Let $\epsilon>0$ and pick $v_1,v_2\in V^c$ so that
 $d_V(x_1,x_2) \leq (1+\epsilon)\log [x_1,x_2;v_1,v_2]$.
To simplify calculations,
we choose a transformation in $GM(\hatH)$ which maps $v_1$
  to zero and $v_2$ to infinity. We recall that this preserves cross-ratios.
By a slight abuse of notation we still write
$x_1,x_2$ for the images in $\hatH$ of the corresponding
points.  We have then $0< d_V(x_1,x_2) \leq (1+\epsilon) \log 
\frac{\|x_2\|}{\|x_1\|}$ so in particular, $\|x_1\|< \|x_2\|$.
When $u_1,u_2\in U$ we have in these
new coordinates,
$\left| \log
\frac{\|u_2\|}{\|u_1\|} \right| = 
\left| \log[u_1,u_2,0,\infty] \right| \leq
 d_V(u_1,u_2) \leq \Delta < +\infty$. In other words,
 $U$ is bounded away from the origin and infinity.

Consider now the formula for
the distance of $x_1,x_2$ relative to $U$. It  splits 
into a sum of two supremums (this splitting is
one of the deeper reasons why the Apollonian metric
is easy to handle):
 \[
    d_U(x_1,x_2) =
         \sup_{u_1\in U^c} \log  \frac{\|x_2-u_1\|}{\|x_1-u_1\|} 
         + \sup_{u_2\in U^c} \log  \frac{\|x_1-u_2\|}{\|x_2-u_2\|}.
\]
The suprema of these two terms are
denoted $\alpha_1$ and $\alpha_2$. They are both finite.
We define the Apollonian ball
  \[
  B_1=B_{\alpha_1}(x_1,x_2)
  = \left\{ u\in \hatH : \mystrut 
        \frac{\|x_1-u\|}{\|x_2-u\|}   < \alpha_1
    \right\} \subset U
  \]
and similarly for the ball $B_2=B_{\alpha_2}(x_2,x_1) \subset U$
 (see Figure \ref{fig cross ratio}).

\begin{figure}
\begin{center}
\epsfig{figure=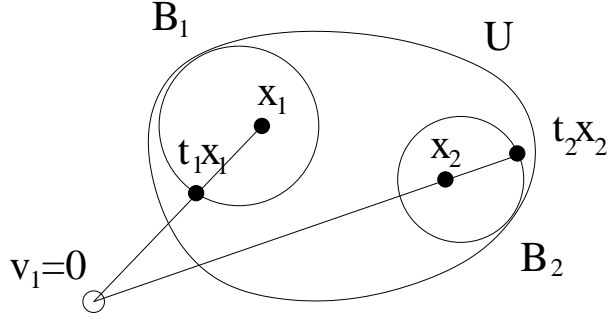,width=8cm}
\end{center}
\caption{Construction of cross-ratios. $v_1=0$ and $v_2=\infty$.}
\mylabel{fig cross ratio}
\end{figure}
A priori $B_1$ is a generalized open ball containing $x_1$
but as $U$ is bounded $B_1$
must be an open ball in the usual bounded sense (and $\alpha_1$
must be greater than one).
Now let $t_1 x_1$
(with $0< t_1< 1$) be the unique intersection of the  segment 
$\{t\;x_1 : 0\leq t\leq 1\}$ and
 the sphere 
$\partial B_{u_1}(x_1,x_2)$. Similarly, let  $t_2 x_2$ 
(with $1< t_2 < +\infty$) be the unique intersection between
the segment 
$\{t\;x_2 : 1\leq t\leq +\infty\}$ and
$\partial B_{u_2}(x_2,x_1)$  (see Figure \ref{fig cross ratio}).
Then 
$t_1 x_1, t_2 x_2 \in \Cl {U}$ and
$\|t_1x_1\|< \|x_1\| < \|x_2\| < \|t_2x_2\|$.  From the
way we defined $t_1$ and $t_2$
we have the following lower bound
\begin{eqnarray}
     d_U(x_1,x_2)  = \alpha_1 + \alpha_2 
    &=& \log \frac{\|x_2-t_1 x_1\|}{\|x_1-t_1 x_1\|} 
           \times \frac{\|x_1-t_2 x_2\|}{\|x_2-t_2 x_2\|}  
    \nonumber \\
    & =& \log  \frac{\|x_2-t_1 x_1\|}
	             {\|x_1\|-\|t_1 x_1\|}
             \times \frac{\|x_1-t_2 x_2\|}
                      {\|t_2 x_2\|-\|x_2\|}
    \nonumber \\
    &\geq&
      \log \frac{\|x_2\|-\|t_1 x_1\|}
              {\|x_1\|-\|t_1 x_1\|}
      \times \frac{\|t_2 x_2\|-\|x_1\|}
                {\|t_2 x_2\|-\|x_2\|}
	      .
      \label{ineq uz}
\end{eqnarray}

The last expression is the cross-ratio of the four (ordered) points
on the positive real line
$0<\|t_1 x_1\| < \|x_1\| < \|x_2\| < \|t_2 x_2\| <+\infty$.
Let us write  $J=(0,\infty)$, 
$K=(\|t_1x_1\|, \|t_2x_2\|)$ and $s_1=\|x_1\|$, $s_2=\|x_2\|$.
By our construction 
$\diam_{J}(K) = \log \frac{\|t_2x_2\|}{\|t_1x_1\|} 
    \leq d_V(t_1x_1,t_2x_2) \leq \diam_V(U) = \Delta$,
where we used that
    $t_1x_1,t_2x_2\in \Cl\; U$
 and that $v_1=0,v_2=\infty\in  V$.
Also $d_{K}(\|x_1\|,\|x_2\|) \leq d_U(x_1,x_2)$
by the above bound (\ref{ineq uz}).
So using Birkhoff's inequality
 (\ref{Birkhoff ineq}) we get
  \[  d_V(x_1,x_2) (1+\epsilon)^{-1} \leq d_{J}(s_1,s_2) \leq 
   \left( \tanh \frac{\diam_{J}(K)}{4} \right) d_{K}(s_1,s_2) \leq
   \left( \tanh \frac{\Delta}{4}\right)  d_U(x_1,x_2),\]
and since $\epsilon>0$ was arbitrary we see that
  \[  d_V(x_1,x_2)\leq 
   \left( \tanh \frac{\Delta}{4}\right)  d_U(x_1,x_2),\]
which is what we aimed to show.\Halmos\\

\section{Some applications}
      \label{sec applic}

In the one dimensional case,
the result of Birkhoff \cite{Bir57} has a vast variety of applications
related to Perron-Frobenius type
of results and the presence of spectral gaps of real operators 
contracting a real convex cone,
see e.g.\ \cite{Bal00}.
In the case of complex operators  similar spectral gap results were
obtained first in \cite{Rug10} and then simplified in \cite{Dub09}
using a complex Hilbert metric and the  2-dimensional version
of the UCP for the Apollonian metric.
We discuss in the following some possible applications in the
case of arbitrary dimension.

\begin{Corollary}
   Let $U \subset V$ and $\Delta$ be as in the Main theorem and write
   $\Gamma(V,U) = \{\gamma \in GM(\hatH): \gamma(V)\subset U\}$
   for the elements of the conformal group that map
    $V$ into $U$. Then for every $\gamma\in \Gamma(V,U)$ we have
    $\gamma^{-1} \in \Gamma(U^c,V^c)$ and the mappings
  $ 
	 \gamma: (V,d_V) \rr (V,d_V)  
	       \  \mbox{and} \ 
	 \gamma^{-1}: ({U^c},d_{U^c}) \rr ({U^c},d_{U^c}) 
      $ 
    are $\left( \tanh \frac{\Delta}{4} \right)$-Lipschitz.
\end{Corollary}
Proof: $\gamma\in \Gamma(U,V)$ preserves cross-ratios,
and $\gamma(V)\subset U$ so writing $\theta=\tanh \Delta/4$
we have for $v_1,v_2\in V$:
 \[ d_V(\gamma(v_1),\gamma(v_2)) \leq 
    \theta \; d_U (\gamma(v_1),\gamma(v_2)) \leq
    \theta \; d_{\gamma(V)} (\gamma(v_1),\gamma(v_2)) \leq
    \theta \; d_{{V}} (v_1,v_2).\]
 The inverse map is bijective so
it maps $U^c$ into $V^c$. We have the same bound for its contraction rate
since
    \[ \diam_V (U) = \diam_{U^c}({V^c}) =
 \; \sup_{v_1,v_2\in V} \;  \sup_{u_1,u_2\in U^c} \;  \log [u_1,u_2;v_1,v_2]. 
\ \ \ \ \Halmos \]
\mbox{}

 \begin{Corollary} In finite dimension when
 $\Cl U \subset \Int V$ for the topology of $(\hatH,\hatd)$,
 then from compactness we see that $\diam_V(U) <+\infty$
 so the embedding $(U,d_U) \hookrightarrow (V,d_V)$ is a
 strict Lipschitz contraction.
 \end{Corollary}

\begin{Lemma}  
  Suppose that
$U\subset B(x_0,R)$,  $R<\infty$.  Then
  \begin{equation}
    \|u_1-u_2\|
    \leq
    \frac{R}{2} \;
    d_U(u_1,u_2)
    , \ \ \forall u_1,u_2\in U.
  \end{equation}
  Suppose that $U\subset V$ and that
  $r=dist(U,V^c)= \sup_{u\in U,w\in V^c} \|u-w\| >0$. Then
  \begin{equation}
    d_V(u_1,u_2)
    \leq \frac{2}{r} \;
    \|u_1-u_2\| 
    , \ \ \forall u_1,u_2\in U.
  \end{equation}
  \end{Lemma}
Proof: 
When $x\in  B(x_0,R)$ and $h$ is small we get 
from a straight-forward calculation:
\[  d_{B}(x,x+h) = \frac{2 R}{R^2 - \|x-x_0\|^2} \|h\| + o(h).\]
Thus, $ds=\frac{2R}{R^2-\|x-x_0\|^2} \|dx\| \geq \frac{2}{R}\|dx\|$ and 
$\|v_1-v_2\| \leq d_B(v_1,v_2) \leq d_V(v_1,v_2)$ (since $V\subset B$).
When $B(u_1,r),B(u_2,r)\subset V$ then for $w\in V^c$:
$\ds \frac{\|u_2-w\|}{\|u_1-w\|} \leq 1 + \frac{\|u_2-u_1\|}{r}$ and
$\ds d_V(u_1,u_2) \leq 2 \log \left( 1 + \frac{\|u_2-u_1\|}{r} \right)
       \leq \frac{2}{r} \|u_2-u_1\|$.
\Halmos\\

 \begin{Theorem}
 Let $U\subset V$ be non-empty proper subsets of 
 $(\hatH,\hatd)$ such that $\Cl V \neq \hatH$ and $\Delta=\diam_V(U)<+\infty$.
 Let $\gamma_1,\ldots,\gamma_k\in \Gamma(V,U)$ and write
 \[ 
       \Lambda\equiv \Lambda(\gamma_1,\ldots,\gamma_k)= 
   \bigcap_{n\geq 1} \Cl \bigcup_{1\leq i_1,\ldots,i_n\leq k} 
      \gamma_{i_1} \circ \cdots \circ \gamma_{i_n} (V)\]
for the associated limit set. 
Then $\Lambda$ is compact and has Hausdorff and Box dimensions not greater
than $-\log k / \log \tanh \frac{\Delta}{4}$.
\end{Theorem}
Proof: Pick $q\in \hatH \setminus \Cl V$ and map $q$ to infinity
by an inversion in $q$.
In the new coordinates $V$ is bounded so by the previous Lemma,
Hilbert distances are bounded by Apollonian distances.
At level $n\geq 1$ each set in the finite union
has diameter not greater than $r=\Delta (\tanh\frac{\Delta}{4})^{n-1}$
which becomes arbitrarily small as $n\rr\infty$.
There are $N_r=k^n$ elements in the union.
As $\Lambda$ is closed and has finite covers of arbitrarily small
diameters it is compact and we have the bound
\[ \dim_H(\Lambda) \leq \limsup_n \ \ \frac{\log N_r}{\log{1/r}}
 \ = \ \frac{\log k } {\log \tanh \frac{\Delta}{4}} \ . \ \ \Halmos \]
 
 When the images $\Cl(\gamma_i(V))$, $1\leq i\leq k$ are
 pairwise disjoint the Hausdorff dimension may also be obtained from
 a Bowen-like formula as in \cite{Rug08} or \cite{MU98}. We omit the details.
 {Note that we do not assume here that $H$ is finite dimensional.}

\begin{Remark}
In finite dimension $d\geq 2$
the Apollonian metric for an open ball $V=B(0,R)$ 
is the same as the hyperbolic metric for the ball,
i.e.\ $ds=2r/(r^2-\|x\|^2) \|dx\|$. In this case it is well-known
that if $\gamma$ maps $V$ inside $V$ and $\gamma(V)$ has
bounded diameter 
then $\gamma$ is a uniform contraction.
\end{Remark}

Other metrics may be constructed from
the Apollonian metric (cf. \cite{Has04}).
Let $v\in H^*$, $\|h\| \leq \|v\|/4$ and write
$x=\frac{\dual{h,v}}{\dual{v,v}} \in [-1/4,1/4]$ and
$\|h\|^2 = \|v\|^2 (x^2 + y^2)$.
Calculus 
 shows that $\left| \frac{1}{2} \log((1+x)^2+y^2) - x\right| \leq
 x^2+y^2$ (when $x\geq -1/4$). Therefore,
\[ 
  \left| \log \frac{\|v+h\|}{\|v\|}-\dual{I(v),h} \right|
  \ =
  \ \left| \log \frac{\|v+h\|}{\|v\|}-\frac{\dual{v,h}}{\dual{v,v}} \right|
    \ \leq \ 
    \frac{\dual{h,h}} {\dual{v,v}}
     .\]
We assume in the following that $U$ is open. 
Let $x\in U$ and set $r=\inf_{u\in U^c} d(x,u)>0$.
When $u_1,u_2\in U^c$ and $\|h\|\leq r/4$ we get:
\[  \left| 
   \log [x,x+h;u_1,u_2] - \dual{I(x-u_1)-I(x-u_2),h} \right|
   \leq 2  \|h\|^2/r^2 .\]
It follows that the following limit exists and define
a Finsler 
(pseudo-) norm on the tangent space of $U$:
 \begin{equation}
      p^\st_{U,x}(h) \equiv \lim_{t\rr 0} \frac{1}{t} d_U(x+th,x) =
          \sup_{u_1,u_2\in U^c}  \left| \dual{ I(x-u_1)-I(x-u_2),h} \right|.
	  \mylabel{pUx}
\end{equation}
 It is only a pseudo-norm
when $U^c$ is contained in a generalized ball, since in
that case $p^\st_{U,x}$ 
may vanish in some directions.
If $\gamma:[0,1] \rr U$ is a continuous path then
we may define its (pseudo-) length to be 
\[ \ell(\gamma)\equiv \limsup_{\delta\rr 0}
\sum_{k=0}^n d_{U}(\gamma(t_{k+1}), \gamma(t_k)),\]
where $0=t_0<t_1<\cdots<t_n=1$ and $t_{k+1}-t_k<\delta$.
Then  
   \begin{equation}
       \dInner_U (x,y) = \inf \{\ell(\gamma): 
       \gamma \in C([0,1],U), \gamma(0)=x,\gamma(1)=y \}
    \end{equation}
defines a (pseudo-)metric which in 
\cite{Has04} was coined the Apollonian inner metric.
   When $\gamma$ is peicewise $C^1$ we have
   $\ell(\gamma)=\int_0^1 p^\st_{U,x} (\dot{\gamma}(t)) \; dt$.
Another possiblity is
to maximize (\ref{pUx}) over directions.
This leads to
a conformal Riemannian
metric $ds=g^\st_U(x) \|dx\|$ with
\begin{equation}
   g^\st_U(x) =
      \sup_{\|h\|=1} p_{U,x}^{\hugestrut}(h) =
      \sup_{u_1,u_2\in U^c} \;\frac{ \|u_1-u_2\|}{\|x-u_1\| \; \|x-u_2\|} .
 \end{equation}
An advantage of this metric
 is perhaps
that it distinguishes points when $U^c$ contains at least two points.
It is easy to see that $g^\st_U(x)$ is continuous (as we assumed
$U$ to be open).
We write $\dRiem_U(x,y)$ for the Riemannian distance of $x$ and $y$
with respect to this metric.

\begin{Corollary}
\mylabel{corol contract} 
    Let $U \subset V \subset \hatH$ 
    (with $\Cl V \neq \hatH$) be non-empty proper subsets
    and let
 $ \Delta =\sup_{u_1,u_2\in U} d_V(u_1,u_2)$
be the diameter of 
the smaller subset within the larger with respect to the Apollonian metric.
Then for every $x,y \in U$:
\begin{equation}
      p^{\strut}_{V,x}(h) \leq \left( \tanh \frac{\Delta}{4} \right)
       \ p^{\strut}_{U,x}(h), \ \ h\in E,
\end{equation}
\begin{equation}
      \dInner_V(x,y) \leq \left( \tanh \frac{\Delta}{4} \right)
       \ \dInner_U(x,y),
\end{equation}
\begin{equation}
      \dRiem_V(x,y) \leq \left( \tanh \frac{\Delta}{4} \right)
      \ \dRiem_U(x,y).
\end{equation}
\end{Corollary}

Proof: 
For $x,x+th\in U$ we have by the Main Theorem
$\frac1t d_V(x,x+th) \leq \tanh\frac{\Delta}{4} \frac1t d_U(x,x+th)$. 
The first inequality
follows. The second follows by taking
limits in the right order.
 For the Riemmannian metric
one has 
      \[
      g^\st_V(x) \leq 
      \sup_{\|h\|=1} p_{U,x}^{\hugestrut}(h) \leq
      \sup_{\|h\|=1} \left( \tanh \frac{\Delta}{4} \right)
      p_{U,x}^{\hugestrut}(h) =
      \left( \tanh \frac{\Delta}{4} \right)
      g^\st_U(x)\]
which yields the last inequality.
      \Halmos\\

\def\theequation{\Alph{section}.\arabic{equation}}
\appendix

\end{document}